\documentclass[aap,MSNbibl,citesort,dvips]{arximspdf}

\doi{10.1214/11-AAP827} 
\volume{22}
\issue{5}
\pubyear{2012}
\firstpage{2086}
\lastpage{2107}

\makeatletter
\newcommand{\eqref}[1]{(\ref{#1})}
\newproclaim{definition}{Definition}
\newtheorem{lemma}[definition]{Lemma}

\newtheorem{theorem}[definition]{Theorem}
\newtheorem{corollary}[definition]{Corollary}
\makeatother

\begin{document}
\begin{frontmatter}

\title{The asymptotic distribution of the length of~Beta-coalescent trees}
\runtitle{Length of Beta-coalescent}

\begin{aug}
\author[A]{\fnms{G\"otz} \snm{Kersting}\corref{}\ead[label=e1]{kersting@math.uni-frankfurt.de}}
\runauthor{G. Kersting}
\affiliation{Goethe Universit\"at}
\address[A]{Goethe Universit\"at Frankfurt\\
Box 187\\
D-60054 Frankfurt/Main\\
Germany\\
\printead{e1}}
\end{aug}

\received{\smonth{8} \syear{2011}}
\revised{\smonth{10} \syear{2011}}

%
\begin{abstract}
We derive the asymptotic distribution of the total length $L_n$ of a
$\operatorname{Beta}(2-\alpha,\alpha)$-coalescent tree for $1<\alpha< 2$, starting
from $n$ individuals. There are two regimes: If $\alpha\le\frac12
(1+ \sqrt5)$, then $L_n$ suitably rescaled has a stable limit
distribution of index $\alpha$. Otherwise $L_n$ just has to be shifted
by a constant (depending on $n$) to get convergence to a nondegenerate
limit distribution. As a consequence, we obtain the limit distribution
of the number $S_n$ of segregation sites. These are points (mutations),
which are placed on the tree's branches according to a Poisson point
process with constant rate.
\end{abstract}

%


\end{frontmatter}
%

\section{Introduction and result}\label{sec1}
In this paper we investigate the asymptotic distribution of the
suitably normalized length $L_n$ of a $n$-coalescent of the
$\operatorname{Beta}(2-\alpha,\alpha)$-type with $1< \alpha< 2$. As a corollary we
obtain the asymptotic distribution of the associated number $S_n$ of
segregating sites, which is the basis of the Watterson estimator \cite
{wa} for the rate $\theta$ of mutation of the DNA. Here we recall that
coalescents with multiple merging such as Beta-coalescents have been
considered in the literature as a model for the genealogical
relationship within certain maritime species~\cite{bobo,elwa}.

Beta-coalescents (and more generally $\Lambda$-coalescents, as
introduced by Pitman~\cite{pi} and Sagitov~\cite{sa}) possess a rich
underlying partition structure, which is nicely presented in detail in
Berestycki~\cite{nber}. For our purposes it is not necessary to recall
all these details, and we refer to the following condensed description
of a $n$-coalescent:

Imagine $n$ particles (blocks in a partition), which coalesce into a
single particle within a random number of steps. This happens in the
manner of a continuous time Markov chain. Namely, if there are
currently $m>1$ particles, then they merge to $l$ particles at a rate
$\rho_{m,l}$ with $1 \le l \le m-1$. Thus
\[
\rho_m= \rho_{m,1} + \cdots+ \rho_{m,m-1}
\]
is the total merging rate, and
\[
P_{m,l} = \frac{\rho_{m,l}}{\rho_m} ,\qquad 1 \le l \le m-1
\]
gives the probability of a jump from $m$ to $l$.

In these models the rates $\rho_{m,l}$ have a specific consistency
structure arising from the merging mechanism. As follows from Pitman
\cite{pi}, they are, in general, of the form
\[
\rho_{m,m-k+1}= \pmatrix{ m\cr{k}} \int_0^1 t^{k-2}(1-t)^{m-k} \Lambda(dt)
,\qquad 2 \le k \le m ,
\]
where $\Lambda(dt)$ is a finite measure on $[0,1]$. The choice $\Lambda
= \delta_0$ corresponds to the original model due to Kingman~\cite{ki},
then $\rho_{m,l}=0$ for $l \ne m-1$. In this paper we assume
\[
\Lambda(dt)= \frac{1}{\Gamma(2-\alpha)\Gamma(\alpha)} t^{1-\alpha
}(1-t)^{\alpha-1} \,dt ,
\]
thus
\[
\rho_{m,m-k+1} = \frac{1}{\Gamma(2-\alpha)\Gamma(\alpha)} \pmatrix{ m\cr
k}
B(k-\alpha,m- k+\alpha) ,
\]
where $B(a,b)$ denotes the ordinary Beta-function. Then the underlying
coalescent is called the $\operatorname{Beta}(2-\alpha,\alpha)$-coalescent. For $\alpha
=1$ it is the Bolthausen--Sznitman coalescent~\cite{bo} and the case
$\alpha\to2$ can be linked with Kingman's coalescent.

The situation can be described as follows: There are the merging times
$0= T_0 < T_1 < \cdots< T_{\tau_n} $
and there is the embedded time discrete Markov chain $n=X_0 > X_1 >
\cdots> X_{\tau_n}=1$, where $X_i$ is the number of particles
(partition blocks) after $i$ merging events, and $\tau_n$ is the number
of all merging events. This Markov chain has transition probabilities
$P_{m,l}$ and, given the event $X_i=m$ with $m>1$, the waiting time
$T_{i+1}-T_i$ to the next jump is exponential with expectation $1/\rho
_{m}$. Since a point process description is convenient later, and we
name the point process
\begin{equation}
\mu_n = \sum_{i=0}^{\tau_n-1} \delta_{X_i}
\label{pp}
\end{equation}
on $\{2,3,\ldots\}$ the \textit{coalescent's point process} downwards
from $n$, abbreviated $\operatorname{CPP}(n)$.

These dynamics can be visualized by a coalescent tree with a root and
$n$ leaves. The leaves are located at height $T_0=0$ and the root at
height $T_{\tau_n}$ above. At height~$T_i$ there are $X_i$ nodes
representing the particles after $i$ coalescing events. The total
branch length of this tree is given by
\begin{equation}
L_n = \sum_{i=0}^{\tau_n-1} X_{i}(T_{i+1}-T_{i}) .
\label{length}
\end{equation}
For $1<\alpha< 2$, the asymptotic magnitude of $L_n$ is obtained by
Berestycki et al. in~\cite{bebesw}; it is proportional to $n^{2-\alpha
}$. The asymptotic distribution of $L_n$ is easily derived for
Kingman's coalescent (see~\cite{de}); it is Gumbel. The case of a
Bolthausen--Sznitman coalescent is treated by Drmota et al.~\cite{dr},
here~$L_n$ properly normalized is asymptotically stable. The case
$0<\alpha< 1$ of a~Beta-coalescent is contained in more general
results of M\"ohle~\cite{mo}. Partial results for the Beta-coalescent
with $1 < \alpha< 2$ have been obtained by Delmas et al.~\cite{de}.

In this paper we derive the asymptotic distribution of the length of
the Beta-coalescent for $1 < \alpha< 2$.
Let $\varsigma$ denote a real-valued stable random variable with index
$ \alpha$, which is normalized by the properties
\begin{equation}
\mathbf{E}(\varsigma)=0 ,\qquad \mathbf{P}( \varsigma> x)=o( x^{-\alpha
}) , \qquad\mathbf{P}( \varsigma<-x)\sim x^{-\alpha} \label{stable}
\end{equation}
for $x \to\infty$. Thus it is maximally skewed among the stable
distributions of index~$\alpha$.

Also let
\[
c_1= \frac{\Gamma(\alpha)\alpha(\alpha- 1)}{2-\alpha} ,\qquad c_2=
\frac{\Gamma(\alpha)\alpha(\alpha-1)^{1+ 1/\alpha}}{\Gamma(2-\alpha
)^{1/\alpha}} .
\]

\begin{theorem}
\label{mainresult}
For the Beta-coalescent with $1< \alpha< 2$:
\begin{longlist}[(iii)]
\item[(i)]
If $1< \alpha< \frac12 (1+ \sqrt5)$ (thus $1+\alpha-\alpha
^2>0$), then
\[
\frac{L_n- c_1n^{2-\alpha}}{n^{1/\alpha+ 1 -\alpha}} \stackrel
{d}{\to} \frac{c_2\varsigma}{(1+\alpha-\alpha^2)^{1/\alpha}} .
\]
\item[(ii)]
If $\alpha=\frac12 (1+ \sqrt5)$, then
\[
\frac{L_n- c_1n^{2-\alpha}}{(\log n)^{1/\alpha} } \stackrel
{d}{\to} c_2\varsigma.
\]
\item[(iii)]
If $\frac12 (1+ \sqrt5) < \alpha< 2$, then
\[
L_n- c_1n^{2-\alpha} \stackrel{d}{\to} \eta,
\]
where $\eta$ is a nondegenerate random variable.
\end{longlist}
\end{theorem}

In fact it is not difficult to see from the proof that $\eta$ has a
density with respect to Lebesgue measure.

This transition at the golden ratio $\frac12 (1+\sqrt5)$ is
manifested in the results of Delmas et al.~\cite{de}. They also show
that the number $\tau_n$ of collisions,\vadjust{\goodbreak} properly rescaled, has an
asymptotically stable distribution. This latter result has been
independently obtained by Gnedin and Yakubovitch~\cite{gn}.

The region within the coalescent tree, where the random fluctuations of
$L_n$ asymptotically arise, are different in the three cases. In case
(i) fluctuations come from everywhere between the root and the leaves,
whereas in case~(iii) they mainly originate at the neighborhood of the
root. Then we have to take care of those summands
$X_{i}(T_{i+1}-T_{i})$ within $L_n$, which have an index~$i$ close to
$\tau_n$. In the intermediate case (ii) the primary contribution stems
from summands with index $i$ such that $\tau_n -n^{1-\varepsilon} \le i
\le\tau_n-n^\varepsilon$ with $0<\varepsilon<\frac12$.

To get hold of these fluctuations, in proving the theorem, we, loosely
speaking, turn around the order of summation in $L_n = \sum_{i=0}^{\tau
_n-1} X_{i}(T_{i+1}-T_{i})$. We shall handle the reversed order by
means of two point processes $\mu$ and~$\nu$ on $\{2,3,\ldots\}$. The
first one, which we call the \textit{coalescent's point process downwards
from} $\infty$ [$\operatorname{CPP}(\infty)$], gives the asymptotic particle numbers
seen from the root of the tree. Here we use Schweinsberg's result \cite
{sw} implying that the Beta-coalescent comes down from infinity for
$1<\alpha< 2$; see~\cite{nber}, Corollary 3.2. (Therefore our method
of proof does not apply to the case of the Bolthausen--Sznitman
coalescent.) The second one is a classical stationary renewal point
process, which can be reversed without difficulty. Two different
couplings establish the links. Thereby the exponential holding times
are left aside at first stage. In this respect, our approach to the
Beta-coalescent differs from others as in Birkner et al.~\cite{bi} or
Berestyki et al.~\cite{be2}. Certainly our proof can be extended to a
larger class of $\Lambda$-coalescents having regular variation with
index $\alpha$ between $1$ and $2$ (compare Definition~4.1 in~\cite
{nber}), which would require some additional technical efforts. It
seems less obvious, whether our concept of a coalescent's point process
downwards from $\infty$ can be realized for a much broader family of
$\Lambda$-coalescents coming down from infinity.

Coalescent trees are used as a model for the genealogical relationship
of $n$ individuals backward to their most recent ancestor. Then one
imagines that mutations are assigned to positions on the tree's
branches in the manner of a~Poisson point process with rate $\theta$.
Let $S_n$ be the number of these segregation sites; see~\cite{nber},
Section 2.3.4. Given $L_n$ the distribution of $S_n$ is Poisson with
mean $\theta L_n$. To get the asymptotic distribution one splits $S_n$
into parts.
\[
S_n - \theta c_1n^{2-\alpha} = (S_n - \theta L_n) + \theta(L_n - c_1
n^{2-\alpha}) .
\]
Since $L_n/c_1 n^{2-\alpha}$ converges to 1 in probability, the first
summand is asymptotically normal and also asymptotically independent
from the second one. Its normalizing constant is $(\theta L_n)^{-1/2} \sim(\theta c_1)^{-1/2}n^{\alpha/2 -1}$. Again there
are two regimes. $n^{1-\alpha/2} = o(n^{1/\alpha+ 1 -\alpha
})$, if and only if $\alpha< \sqrt2$. Partial results are contained
in Delmas et al.~\cite{de}. We obtain

\begin{corollary} Let $\zeta$ denote a standard normal random variable,
which is independent of $\varsigma$.
\begin{longlist}[(iii)]
\item[(i)]
If $1< \alpha< \sqrt2$, then
\[
\frac{S_n- \theta c_1n^{2-\alpha}}{n^{1/\alpha+ 1 -\alpha}}
\stackrel{d}{\to} \frac{\theta c_2\varsigma}{(1+\alpha-\alpha
^2)^{1/\alpha}} .
\]
\item[(ii)]
If $\alpha=\sqrt2$, then
\[
\frac{S_n- \theta c_1n^{2-\alpha}}{n^{1- \alpha/2}} \stackrel
{d}{\to} \sqrt{\theta c_1} \zeta+\frac{\theta c_2\varsigma}{(1+\alpha
-\alpha^2)^{1/\alpha}} .
\]
\item[(iii)]
If $\sqrt2 < \alpha< 2$, then
\[
\frac{S_n- \theta c_1n^{2-\alpha}}{n^{1- \alpha/2}} \stackrel
{d}{\to} \sqrt{\theta c_1} \zeta.
\]
\end{longlist}
\end{corollary}

This is the organization of the paper: Section~\ref{sec2} contains an elementary
coupling of two $\mathbb{N}$-valued random variables. It is used in
Section~\ref{sec3}, where we introduce and analyze coalescent's point processes,
and in Section~\ref{sec4}, where we couple these point processes to stationary
point processes. Section~\ref{sec5} assembles two auxiliary results on sums of
independent random variables. Finally the proof of Theorem \ref
{mainresult} is given in Section~\ref{sec6}.

\section{A coupling}\label{sec2}

In this section, let the natural number $m$ be fixed. We introduce a
coupling of the transition probabilities $P_{m,l}$ and a distribution,
which does not depend on $m$.
From the representation of the Beta-function by means of the $\Gamma
$-function and its functional equation, we have
\begin{eqnarray*}
\rho_{m,m-k+1} &=& \frac1{\Gamma(2-\alpha) \Gamma(\alpha
)}\frac{m!}{\Gamma(m)} \frac{\Gamma(k-\alpha)}{k!} \frac{\Gamma
(m-k+\alpha)}{(m-k)!} \\
& =& \frac1{\Gamma(2-\alpha) \Gamma(\alpha)} \frac{\Gamma(k-\alpha
)}{\Gamma(k+1)} \frac{(m-k+1)\cdots m}{(m-k+\alpha)\cdots(m-1+\alpha
)} \frac{\Gamma(m+\alpha)}{\Gamma(m)} ,
\end{eqnarray*}
thus
\[
P_{m,m-k} = d_{mk} \frac{\Gamma(k+1-\alpha)}{\Gamma(k+2)} ,\qquad k
\ge1
\]
with
\[
d_{mk}=d_m \frac{(m-k )\cdots(m-1)}{(m+\alpha-k-1)\cdots(m +\alpha
-2)}
\]
and a normalizing constant $d_m>0$ (also dependent on $\alpha$). Recall
from the \hyperref[sec1]{Introduction} that given $X_0=m$ the quantities $P_{m,m-k}$ are
the weights of the distribution of the downward jump $U=X_0-X_1$. For a
more detailed discussion of this ``law of first jump,'' we refer to
Delmas et al.~\cite{de}.

It is natural to relate this distribution to the distribution of some
random variable~$V$ with values in $\mathbb{N}$ and distribution given by
\begin{equation}
\mathbf{P}(V=k) = \frac{\alpha}{\Gamma(2-\alpha)} \frac{\Gamma
(k+1-\alpha)}{\Gamma(k+2)} ,\qquad k \ge1 . \label{Z1}\vadjust{\goodbreak}
\end{equation}
This kind of distribution appears for Beta-coalescents already in
Bertoin and Le Gall~\cite{bega} (see their Lemma 4), in Berestycki et
al.~\cite{be2} (in the context of frequency spectra) as well as in
Delmas et al.~\cite{de}.
There the normalizing constant is determined and the following formulas derived:
\begin{equation}
\mathbf{E}(V) =\frac1{ \alpha-1}\quad \mbox{and}\quad \mathbf{P}(V
\ge k)= \frac1{\Gamma(2-\alpha)} \frac{\Gamma(k+1-\alpha)}{\Gamma
(k+1)} . \label{Z2}
\end{equation}
From Stirling's approximation,
\begin{equation}
\mathbf{P}(V=k) \sim\frac{\alpha}{\Gamma(2-\alpha)} k^{-\alpha-1}
\quad\mbox{and}\quad \mathbf{P}(V\ge k) \sim\frac1{\Gamma(2-\alpha
)} k^{-\alpha} . \label{Z3}
\end{equation}

The sequence $d_{mk}$ is decreasing in $k$ for fixed $m$, and thus the
same is true for $P_{m,m-k}/\mathbf{P}(V=k)$. Therefore $V$
stochastically dominates the jump size $U$, that is, for all $k\ge1$,
\begin{equation}
\mathbf{P}(U\ge k \mid X_0=m) \le\mathbf{P}(V \ge k) .
\label{dominated}
\end{equation}
We like to investigate a coupling of $U$ and $V$, where $U \le V$ a.s.
It is fairly obvious that this can be achieved in such a way that
\begin{equation}
\mathbf{P}(U= k \mid V=k)= 1 \wedge\frac{P_{m,m-k}}{\mathbf{P}(V=k)} =
1 \wedge\frac{d_{mk}}d .
\label{couple}
\end{equation}
[Indeed one may put
\[
\mathbf{P}(U=j \mid V=k) = \biggl( 1- \frac{P_{m,m-k}}{\mathbf{P}(V=k)}
\biggr)^+ \frac{(P_{m,m-j}-\mathbf{P}(V=j))^+}{\mathbf{P}(U<k_m)-\mathbf{P}(V<k_m)}
\]
for $j \neq k$ with $k_m= \min\{ k \ge1\dvtx P_{m,m-k}\le\mathbf{P}(V=k)\}
$. There are other possibilities; later it will only be important that
we commit to one of them.]

\begin{lemma}\label{fehlerw}
For a coupling $(U,V)$ fulfilling~\eqref{couple}, it holds
\[
\mathbf{P}(U\neq V) \le\frac1{(\alpha-1)m}
\quad\mbox{and}\quad
\mathbf{P}(V\ge k \mid U\neq V) \le c k^{1-\alpha}
\]
for all $k\ge1$ and some $c < \infty$, which does not depend on $m$.
\end{lemma}

\begin{pf} Because of $\alpha< 2$,
\[
\frac{(m-k )\cdots(m-1)}{(m+\alpha-k-1)\cdots(m +\alpha-2)} \ge\frac
{(m-k )\cdots(m-1)}{(m-k+1)\cdots m}= \frac{m-k}m,
\]
and because of $\alpha>1$
\[
\frac{(m+\alpha-k-1)\cdots(m +\alpha-2)}{(m-k )\cdots(m-1)} \ge
\biggl(\frac{m+\alpha-1}m\biggr)^k \ge1+ k\frac{\alpha-1}m ,
\]
consequently
\[
{1-\frac km}\le\frac{d_{mk}}{d_m} \le\frac1{1+ (\alpha-1)
k/m} .
\]
It follows
\[
\biggl( 1-\frac km\biggr) \mathbf{P}(V=k) \le\frac d {d_m}P_{m,m-k} \le
\mathbf{P}(V=k)
\]
for all $k \ge1$ with $d= \alpha/\Gamma(2-\alpha)$. Summing over $k$ yields
\[
1 - \frac1{m} \mathbf{E}(V) \le\frac d{d_m} \le1 \quad\mbox{or}\quad
1 \le\frac{d_m}d \le\frac1{(1- 1/{((\alpha-1)m)})^+} .
\]
Combining the estimates, we end up with
\begin{equation}
1-\frac km \le\frac{d_{mk}} {d}\le\frac1{(1+ (\alpha-1)
k/m)(1- 1/{((\alpha-1)m)})^+}
\label{couple2}
\end{equation}
for all $k \ge1$.

Now from~\eqref{couple},~\eqref{couple2}
\begin{eqnarray*}
\mathbf{P}(U\neq V) &=& \sum_{k \ge1}\bigl (\mathbf{P}(V=k)-P_{m,m-k}\bigr)^+ \\
 &=&
\sum_{k\ge1} \mathbf{P}(V=k) \biggl(1- \frac{d_{mk}}d\biggr)^+ \le\sum_{k \ge1}
\mathbf{P}(V=k) \frac km ,
\end{eqnarray*}
and thus from~\eqref{Z2},
\[
\mathbf{P}(U\neq V) \le\frac1{(\alpha- 1)m}
\]
which is our first claim.

Also, letting $ m \ge2/(\alpha-1)$ and $k' =2\lceil(\alpha
-1)^{-2}+(\alpha-1)^{-1}\rceil$, then
\begin{eqnarray*}
\biggl(1+ (\alpha-1) \frac{k'}m\biggr)\biggl(1- \frac1{(\alpha
-1)m}\biggr)^+&=& 1+ (\alpha-1)\frac{k'}m - \frac1{(\alpha-1)m}- \frac
{k'}{m^2}\\
&\ge&1+ \frac{\alpha-1}2 \frac{k'}m - \frac1{(\alpha-1)m} \ge1+
\frac1m .
\end{eqnarray*}
From~\eqref{couple2},
\[
1 - \frac{d_{mk'}}d \ge1- \frac1{1+1/m} \ge\frac1{2m},
\]
and from~\eqref{couple},
\[
\mathbf{P}(U\neq V) \ge\mathbf{P}(U\neq k',V=k')= \biggl(1- \frac
{d_{mk'}}d\biggr)^+\mathbf{P}(V=k') \ge\frac1{2m} \mathbf{P}(V=k')
\]
for $m \ge2/(\alpha-1)$. It follows that there is a $\eta>0$ such
that for all $m \ge1$
\[
\mathbf{P}(U\neq V) \ge\frac1{\eta m } .
\]
Now from~\eqref{couple} and~\eqref{couple2},
\begin{eqnarray*}
\mathbf{P}(V=k \mid U\neq V) &=& \frac{\mathbf{P}( U\neq
k,V=k)}{\mathbf{P}(U\neq V)} \\
 &= &\frac{(1- {d_{mk}}/d)^+}{\mathbf{P}(U\neq V)}\mathbf{P}(V=k) \le\eta k \mathbf{P}(V=k),
\end{eqnarray*}
and the second claim follows from~\eqref{Z3}.
\end{pf}

\section{The coalescent's point process}\label{sec3}

Let $\mu$ denote a point process on $\{2,3, \ldots\}$. For any interval
$I$, let $\mu_I$ be the point process on $\{2,3, \ldots\}$ given by
\[
\mu_I(B)= \mu(B \cap I) ,\qquad B \subset\{2,3, \ldots\} .
\]

We call $\mu$ a \textit{coalescent's point process} downwards from
$\infty
$, shortly a~$\operatorname{CPP}(\infty)$, if the following properties hold:
\begin{itemize}
\item
$\mu(\{2,3,\ldots\})=\infty$ and $\mu(\{n\})=0$ or $1$ for any $n \ge
2$ a.s.
\item
For $n\ge2$ we have that, given the event $\mu(\{n\})=1$, and given $\mu
_{[n+1,\infty)}$ the point process $\mu_{[2,n]}$ is a $\operatorname{CPP}(n)$ a.s.
\end{itemize}
Recall that a point process is called a $\operatorname{CPP}(n)$, if it can be
represented as in~\eqref{pp}.

\begin{theorem}\label{cpp}
Let $1<\alpha< 2$. Then the $\operatorname{CPP}(\infty)$ exists and is
unique in distribution.
\end{theorem}

We prepare the proof by two lemmas.

\begin{lemma} \label{le1}
Let $\mu$ be a $\operatorname{CPP}(n)$ with $1 < n \le\infty$. Then for any
$\varepsilon>0$ there is a natural number $r$ such that for any
interval $I=[a,b] $ with $2 \le a < b <n$ and $b-a \ge r$, we have
\[
\mathbf{P}\bigl(\mu(I) =0\bigr) \le\varepsilon.
\]
\end{lemma}

\begin{pf}
For $I=[a,b]$,
\[
\{ \mu(I)=0\} = \bigcup_{m=b+1}^n \bigl\{\mu(\{m\})=1, \mu([a,m-1])=0
\bigr\}\qquad \mbox{a.s.},
\]
since $\mu(\{n\})=1$ for $n<\infty$ and $\mu(\{2,3,\ldots\})=\infty$
a.s. for $n=\infty$. Thus from~\eqref{pp},
\[
\mathbf{P}\bigl(\mu(I) =0\bigr) \le\sum_{m=b+1}^n \mathbf{P}(X_1 < a \mid X_0=m) .\vadjust{\goodbreak}
\]
Applying~\eqref{dominated} to $U=X_0-X_1$ it follows that
\begin{equation} \label{jumpsize}
\mathbf{P}\bigl(\mu(I) =0\bigr) \le\sum_{m=b+1}^{n } \mathbf{P}(V> m-a) \le\sum
_{k=1}^\infty\mathbf{P}(V>b-a+k) .
\end{equation}
Since $\mathbf{E}(V)<\infty$, this series is convergent and the claim follows.
\end{pf}

The next lemma prepares a coupling of CPPs.

\begin{lemma} \label{le2}
Let $\mu,\mu'$ be two independent CPPs coming down from $n,n' \le\infty
$. Then for any $\varepsilon>0$ there is a natural number $s$ such
that for any $b$ sufficiently large and $n,n' > b$, we have
\[
\mathbf{P}\bigl( \mu(\{j\})=\mu'(\{j\})=1 \mbox{ for some } j =b-s,\ldots
,b \bigr) \ge1-\varepsilon.
\]
\end{lemma}

\begin{pf}
First let $n<\infty$. We construct a coupling of a $\operatorname{CPP}(n)$ $\mu$ to an
i.i.d. random sequence. Consider random variables $U_1,V_1,U_2,V_2,
\ldots$ and $n=X_0,X_1, \ldots$ with $X_i=n-U_1- \cdots-U_i$, which
are constructed inductively as follows: If $U_1,V_1, \ldots, U_i,V_i$
are already gotten, then given the values of these random variables let
$V_{i+1}$ be a copy of the random variable $V$ from Section~\ref{sec2} and
couple $U_{i+1}$ to $V_{i+1}$ as in Section~\ref{sec2}, with $m=X_i$. For
definiteness, put $U_{i+1}=0$ if $X_i=1$. Then $V_1,V_2, \ldots$ are
i.i.d. random variables with distribution~\eqref{Z1}, and $X_0>X_1 >
\cdots> X_{\tau_n-1}$ are the points of a $\operatorname{CPP}(n)$ $\mu$ down from
$n$, where $\tau_n$ is the natural number $i$ such that $X_i=1$ for the
first time.

Now let $k$ be a natural number. Then $X_{i-1} \ge n- U_1-\cdots-U_k$
for $i\le k$.
Thus for any $\eta>0$ and $n \ge6k\eta^{-1}\mathbf{E}(V) +2$, from
Lemma~\ref{fehlerw},
\begin{eqnarray*}
&&\mathbf{P}\bigl( U_i\neq V_i \mbox{ for some } i\le k, U_1+ \cdots+ U_k
\le6k \eta^{-1}\mathbf{E}(V)\bigr)\\
&&\qquad \le\sum_{i=1}^k \mathbf{P}\bigl(U_i\neq V_i, X_{i-1}\ge n- 6k\eta
^{-1}\mathbf{E}(V)\bigr)\le\frac k{(\alpha-1)(n -6k\eta^{-1}\mathbf{E}(V))} ,
\end{eqnarray*}
thus
\[
\mathbf{P}\bigl(U_i\neq V_i \mbox{ for some } i\le k, U_1+ \cdots+ U_k
\le6k \eta^{-1}\mathbf{E}(V)\bigr) \le\frac\eta6
\]
if $n$ is large enough. Also $\mathbf{E}(U_i) \le\mathbf{E}(V)$ because
of~\eqref{dominated}. Thus from Markov's inequality,
\begin{equation}
\mathbf{P}\bigl(U_1+\cdots+U_k > 6k\eta^{-1}\mathbf{E}(V)\bigr) \le\frac\eta6 ,
\label{markov}
\end{equation}
and consequently
\[
\mathbf{P}(U_i\neq V_i \mbox{ for some } i\le k) \le\frac\eta3
\]
if $n$ is sufficiently large (depending on $\eta$ and $k$).\vadjust{\goodbreak}

Next let $l$ be a natural number and $n'=n+l$. Let
$U_1',V_1',U_2',V_2', \ldots$ and $n'=X_0',X_1', \ldots$ an analog
construction with random variables, which are independent of
$U_1,V_1,U_2,V_2, \ldots.$ Then also
\[
\mathbf{P}(U_i'\neq V_i' \mbox{ for some } i\le k) \le\frac\eta3 .
\]
Moreover because $V$ has finite expectation and because of independence
from classical results on recurrent random walks,
\[
\mathbf{P}\Biggl(\sum_{i=1}^jV_i \neq\sum_{i=1}^j V_i'-l \mbox{ for all
}j \le k\Biggr) \le\frac\eta6 ,
\]
if only $k$ is sufficiently large (depending on $l$). Combining the
estimates we obtain
\[
\mathbf{P}\Biggl(\sum_{i=1}^j U_i \neq\sum_{i=1}^j U_i'-l \mbox{ for all
}j \le k\Biggr) \le\frac{5\eta}6 .
\]
For the corresponding independent CPPs $\mu$ and $\mu'$ coming down
from $n$ and $n'=n+l$, this implies, together with~\eqref{markov},
\[
\mathbf{P}\bigl( \mu(\{j\})= \mu'(\{j\})=1 \mbox{ for some } j \in
[n-6k\eta^{-1}\mathbf{E}(V),n]\bigr) \ge1- \eta.
\]

Leaving aside the coupling procedure we have proved the following: Let
$\eta> 0$, let $l$ be a natural number and let $\mu$ and $\mu'$
denote independent CPPs coming down from $n<\infty$ and $n'=n+l$. Then
there is a natural number $r'$ such that
\begin{equation}
\mathbf{P}\bigl(\mu(\{j\})= \mu'(\{j\})=1 \mbox{ for some } j=n-r',
\ldots,n \bigr) \ge1 - \eta,
\label{forsome}
\end{equation}
if only $n$ is large enough.

With this preparation we come to the proof of the lemma. Let
$\varepsilon>0$, $b \ge2$ and let $ n,n'>b$. Denote
\[
M= \max\bigl\{k \le b\dvtx \mu(\{k\})=1 \bigr\} , \qquad M'= \max\bigl\{k \le
b\dvtx \mu'(\{k\})=1 \bigr\}
\]
(with the convention $M=1$, if $\mu([2,b])=0$). From Lemma~\ref{le1},
\[
\mathbf{P}( M,M' \in[b-r,b]) \ge1-\frac\varepsilon2
\]
for some $r$ and $b>r+2$. Then
\begin{eqnarray*}
&&\mathbf{P}\bigl(\mu(\{j\})= \mu'(\{j\})=1 \mbox{ for no } j\in[b-r'-r,b]
\bigr) \\
&&\qquad\le\frac\varepsilon2+ \mathbf{P}\bigl(\mu(\{j\})= \mu'(\{j\}
)=1 \\
&&\hspace*{42pt}\qquad\mbox{for no } j\in[b-r'-r,b] ;b-r \le M,M' \le b\bigr) \\
&&\qquad\le\frac\varepsilon2+ 2\sum_{b-r\le m< m'\le b} \mathbf{P}\bigl(\mu
(\{j\})= \mu'(\{j\})=1 \\
&&\hspace*{105pt}\qquad \mbox{for no } j=m-r', \ldots,m \mid X_0=m,X_0'=m'\bigr).
\end{eqnarray*}
From~\eqref{forsome} it follows that the right-hand probabilities are
bounded by $\eta= \varepsilon/4r^2$, if~$b$ is only sufficiently large. Then
\[
\mathbf{P}\bigl(\mu(\{j\})= \mu'(\{j\})=1 \mbox{ for no } j\in[b-r'-r,b]
\bigr) \le\varepsilon,
\]
which is our claim with $s=r+r'$.
\end{pf}

As a corollary, we note:

\begin{lemma}
Let $\mu$ and $\mu'$ be two independent $\operatorname{CPP}(\infty)$. Then a.s. $\mu(\{
j\})=\mu'(\{j\})=1$ for infinitely many $j \in\mathbb{N}$.
\end{lemma}

\begin{pf}
From the preceding lemma there are numbers $b_1<b_2< \cdots$ such that
\[
\mathbf{P}\bigl( \mu(\{j\})=\mu'(\{j\})=1 \mbox{ for no } j =b_k,\ldots
,b_{k+1} \bigr) \le2^{-k} .
\]
Now an application of the Borel--Cantelli lemma gives the claim.
\end{pf}

\begin{pf*}{Proof of Theorem \protect\ref{cpp}}
The existence follows from the fact that for $ \alpha>1$ the
corresponding Beta-coalescent $(\Pi_t)_{t \ge0}$ comes down from
infinity~\cite{sw}, which means that the number of blocks in $\Pi_t$ is
a finite number $N_t$ for each $t>0$. Put $\mu(\{k\})=1$, if and only
if $N_t=k$ for some $t>0$.

Uniqueness follows from the last lemma and a standard coupling argument.
\end{pf*}

\section{A bigger coupling}\label{sec4}

Now let $\nu$ be a stationary renewal point process on $\{2,3,\ldots\}
$; that is, if we denote the points of $\nu$ by $ 2 \le R_1 < R_2 <
\cdots$, then the increments $R_{i+1}-R_i$ are independent for $i \ge
0$ (with $R_0=1$) and $R_{i+1}-R_i$ has for $i\ge1$ the distribution
\eqref{Z1}. A stationary version of the process exists, since $\mathbf{E}(V)<\infty$, such that the distribution of $R_1$ may be adjusted in
the usual way to obtain stationarity, that is,
\begin{equation}
\mathbf{P}(R_1= r) = \frac{\mathbf{P}(V\ge r-1)}{\mathbf{E}(V)} ,\qquad
r=2,3,\ldots.
\label{stat}
\end{equation}
Stationarity is of advantage for us. Then $\nu$ may be considered as
restriction of a stationary point process on $\mathbb{Z}$. Such a
process is invariant in distribution under the transformation $z
\mapsto z_0-z$, $z \in\mathbb{Z}$ with $z_0 \in\mathbb{Z}$. Therefore
$\nu$, restricted to $\{2,\ldots,n\}$ looks the same, when considered
upwards or downwards.

In this section we introduce a coupling between $\nu$ and the
$\operatorname{CPP}(\infty)$ $\mu$, which allows us later to replace $\mu$ by $\nu$.
Given $b \ge2$ let, as above,
\[
M= \max\bigl\{k \le b \dvtx \mu(\{k\})=1\bigr\} ,\qquad M' =\max\bigl\{k \le b \dvtx \nu(\{
k\})=1\bigr\} .
\]
Again, if there is no $k\le b$ such that $\mu(\{k\})=1$, we put $M =1$,
and similary for $M'$.
Let $\lambda_b$ and $\lambda_b'$ denote the distributions of $ M $ and
$ M'$ (both dependent on $b$).\vadjust{\goodbreak}

Now for $r \in\mathbb{N}$ we consider the following construction of $\mu
$ and $\nu$, restricted to $[2^{r-1}+1,2^r]$. Take any coupling
$(M,M')$ of $\lambda_{2^r}$ and $\lambda_{2^r}'$. Given $(M,M')$
construct random variables $U_1,V_1,U_2,V_2, \ldots$ inductively as in
the proof of Lemma~\ref{le2}, using the coupling of Section~\ref{sec2}. Here we
start with $X_0=M$. Also let $Y_0=M'$,
\begin{equation} \label{bigcoupl1}
X_i= M-U_1- \cdots- U_i ,\qquad Y_i= M'-V_1- \cdots- V_i ,\qquad
i \ge1 ,
\end{equation}
and
\begin{equation} \label{bigcoupl1a}
N= \min\{i\ge0 \dvtx X_i \le2^{r-1}\} ,\qquad N' = \min\{ i \ge0 \dvtx Y_i
\le2^{r-1} \} .
\end{equation}
The whole construction is interrupted at the moment $N\vee N'$. Maybe
$M,M' \le2^{r-1}$, then no step of the construction is required.
Clearly the following statements are true:
\begin{itemize}
\item
The point process $\sum_{i=0}^{N-1} \delta_{X_i}$ is equal in
distribution to $\mu$, restricted to $[2^{r-1}+1,2^r]$.
\item
The point process $\sum_{i=0}^{N'-1} \delta_{Y_i}$ is equal in
distribution to $\nu$, restricted to $[2^{r-1}+1,2^r]$.
\item
$X_N$ and $Y_{N'}\vee1$ have the distributions $\lambda_{2^{r-1}}$ and
$\lambda_{2^{r-1}}'$.
\end{itemize}
The complete coupling is
\begin{eqnarray}\label{bigcoupl2}
\Phi^r(M,M')&=& \Biggl(\sum_{i=0}^{N-1} \delta_{X_i},\sum_{i=0}^{N'-1}
\delta_{Y_i},X_N,Y_{N'}\vee1\Biggr)
\nonumber
\\[-8pt]
\\[-8pt]
\nonumber
&=& (\phi^r_1,\phi^r_2,\phi^r_3,\phi^r_4)\qquad \mbox{(say)} .
\end{eqnarray}
Its distribution is uniquely determined by the distribution of the
coupling $(U,V)$ from Section~\ref{sec2}. The following continuity property is obvious:
\begin{itemize}
\item If we have a sequence $(M_n,M_n')$ of couplings of $\lambda
_{2^r}$ and $\lambda_{2^r}'$ such that $ (M_n,M_n') \stackrel{d}{\to}
(M,M')$, then $(M,M')$ is also a coupling of $\lambda_{2^r}$ and
$\lambda_{2^r}'$ and
\[
\Phi^r(M_n,M_n') \stackrel{d}{\to} \Phi^r(M,M') .
\]
\end{itemize}

Another obvious fact is that this construction can be iterated: Given
$\Phi^r(M,M')$ we construct $\Phi^{r-1}(\phi^r_3,\phi^r_4)$ and so
forth. Thus starting with the independent coupling $(M,M')$ (i.e., $M$
and $M'$ are independent), we obtain the tupel
\[
\Psi^r= (\Phi^{1,r}(M_{1,r},M_{1,r}'), \Phi
^{2,r}(M_{2,r},M_{2,r}'), \ldots, \Phi^{r,r}(M_{r,r},M_{r,r}')) ,
\]
where $(M_{r,r},M_{r,r}')=(M,M')$ and $(M_{s,r},M_{s,r}')=(\phi^{
s+1,r}_3,\phi^{s+1,r}_4)$ for $s<r$. Since $M_{s,r}$ and $M_{s,r}'$ are
no longer independent in general, the tupels $\Psi^r$ are initially not
consistent for different $r$. To enforce consistency note that for
fixed $s$ the distributions of $(M_{s,r},M_{s,r}')$ are tight for $r
\ge s$, since they take values in the finite set\vadjust{\goodbreak} $\{1,\ldots,2^s\}\times
\{1,\ldots,2^s\}$. Thus by a diagonalization argument, we may obtain a
sequence $1\le r_1<r_2 < \cdots$ such that
\[
(M_{s,r_n},M_{s,r_n}) \stackrel{d}{\to} (M_{s,\infty},M_{s,\infty}')
\]
for certain couplings $(M_{s,\infty},M_{s,\infty}')$ of $\lambda_{2^s}$
and $\lambda_{2^s}'$.

If we make use instead of the independent coupling $(M,M')$, now
$(M_{r,\infty},\break M_{r,\infty}')$ as starting configuration in the
construction of $\Psi^r$, then we gain consistency in the sense that
\[
\Psi^{r-1} \stackrel{d}{=} (\Phi^{1,r}(M_{1,r},M_{1,r}'), \Phi
^{2,r}(M_{2,r},M_{2,r}'), \ldots, \Phi^{r-1,r}(M_{r-1,r},M_{r-1,r}')) .
\]
Proceeding to the projective limit, we obtain the ``big coupling,''
\begin{equation} \Psi^{\infty} =(\Phi^{1,\infty}(M_{1, \infty},M_{1,
\infty}'), \Phi^{2, \infty}(M_{2, \infty},M_{2, \infty}'), \ldots)
.
\label{bigcoupl6}
\end{equation}
It has the property that
\begin{equation}
\mu= \sum_{r=1}^\infty\phi_1^{r,\infty} \quad\mbox{and}\quad \nu= \sum
_{r=1}^\infty\phi_2^{r,\infty}
\label{bigcoupl7}
\end{equation}
are coupled copies of our $\operatorname{CPP}(\infty)$ and stationary point process.

In order to estimate the difference between both point processes, we go
back to~\eqref{bigcoupl1},~\eqref{bigcoupl1a} and estimate the tail of the
distribution of

\begin{equation}
D_r= \max_{i \le N \wedge N'} |X_i-Y_i| .
\label{bigcoupl3}
\end{equation}

\begin{lemma} \label{bigcoupl4}
There is a constant $c>0$ such that for all $r \ge1$ and all $t >0$,
\[
\mathbf{P}( D_r >t) \le c t^{1-\alpha}.
\]
\end{lemma}

\begin{pf}
For $i \le N\wedge N'$, we have
\begin{eqnarray}\label{Dr0}
|X_i-Y_i| &\le&\sum_{j \le N\wedge N'} |U_j-V_j| + |X_0-Y_0|
\nonumber
\\[-8pt]
\\[-8pt]
\nonumber
&\le&\sum_{j \le N\wedge N'} |U_j-V_j| + (2^r-M)+(2^r-M') .
\end{eqnarray}
From~\eqref{Z3},~\eqref{jumpsize},
\begin{equation}
\mathbf{P}(2^r-M > t) \le\sum_{k \ge t} \mathbf{P}( V \ge k) \le
ct^{1-\alpha}
\label{Dr1}
\end{equation}
for a suitable $c>0$.

Because of stationarity $2^r-M'$ and $(R_1-2)\wedge(2^r-1)$ are equal
in distribution, therefore because of~\eqref{Z3},~\eqref{stat}
\begin{equation}
\mathbf{P}(2^r- M' >t ) \le\mathbf{P} (R_1 >t) \le ct^{1-\alpha}
\label{Dr2}
\end{equation}
for a suitable $c>0$.

Finally from Lemma~\ref{fehlerw} $U_j\neq V_j$ occurs for $j \le N$ at
most with probability $p=2^{1-r}/(\alpha-1)$ and then\vadjust{\goodbreak} $|U_j-V_j| \le
V_j$ a.s. Also because of Lemma~\ref{fehlerw} these~$V_j$ can be
stochastically dominated by random variables $a+ b\zeta_j$ with
constants $a,b> 0$ and positive i.i.d. random variables $\zeta_j$,
which possess a~stable distribution of index $\alpha- 1$ and Laplace
transform $\exp(- \lambda^{\alpha-1})$. Also $N\wedge N' \le2^{r-1}=w$
(say). Thus $\sum_{j \le N\wedge N'} |U_j-V_j|$ is stochastically
dominated by the random variable
\[
W=\sum_{j=0}^{w} (a+b \zeta_j) I_{j},
\]
where $I_j$ are i.i.d. Bernoulli with success probability $p$.
Let $\varphi(\lambda)= \exp(-a\lambda-(b \lambda)^{\alpha-1})$ be the
Laplace transform of $a+b\zeta_j$. Then $W$ has the Laplace transform
\[
\sigma(\lambda) = \bigl(1-p\bigl(1- \varphi(\lambda)\bigr) \bigr)^w.
\]
It follows $1- \sigma(\lambda) \le wp (1-\varphi(\lambda)) \le
(1-\varphi(\lambda)) /(\alpha- 1)$. From the well-known identity
$\lambda\int_0^\infty e^{-\lambda x} \mathbf{P}(W>x) \,dx = 1-\sigma
(\lambda) $, it follows that
\begin{eqnarray*}
e^{-1} \mathbf{P}(W >t) &\le& t^{-1}\int_0^\infty e^{-x/t} \mathbf{P}(W>x)
\,dx \\[-2pt]
&=& 1-\sigma(1/t) \le\frac1{\alpha-1}\bigl(1-\exp
\bigl(-at^{-1}-(bt)^{1-\alpha}\bigr)\bigr) .
\end{eqnarray*}
Thus
\begin{equation}
\mathbf{P}\biggl(\sum_{j \le N\wedge N'} |U_j-V_j|>t \biggr) \le\mathbf{P}(W >t) \le ct^{1-\alpha}
\label{Dr3}
\end{equation}
for a suitable $c>0$.
Using estimates~\eqref{Dr1} to~\eqref{Dr3} in~\eqref{Dr0} yields our
claim.\vspace*{-3pt}~%
\end{pf}

Additionally, we note that
\begin{equation}
|N-N'| \le D_r .
\label{bigcoupl5}
\end{equation}
Indeed, if $N<N'$, then $X_N\le2^{r-1} $, thus $Y_N\le2^{r-1}+D_r$.
Further $Y_{N'-1} > 2^{r-1}$, which implies $N'-1-N\le Y_N-Y_{N'-1} \le
D_r-1$. The case $N'<N$ is treated in the same way.\vspace*{-2pt}

\section{On sums of independent random variables}\label{sec5}

The following lemma can be deduced from well-known results (see, e.g.,
Petrov~\cite{pe}), but a direct proof seems more convenient. Let
\[
\gamma= \frac1{\alpha- 1} .\vspace*{-2pt}
\]

\begin{lemma} \label{le51}
Let $V_1,V_2, \ldots$ be i.i.d. copies of the random variable \eqref
{Z1}. Then for any $\beta\in\mathbb R$ and any $\varepsilon> 0$ a.s.,
\[
\sum_{k=1}^n k^{-\beta} (V_k-\gamma) = \eta_n+ o(n^{1/\alpha-
\beta+ \varepsilon}),
\]
where $\eta_n$ is a.s. convergent.\vadjust{\goodbreak}
\end{lemma}

\begin{pf} Let $\varepsilon>0$.
A short calculation gives that
$\mathbf{E}(V_k^2 ; V_k \le k^{1/\alpha+ \varepsilon})$ is of
order $k^{2/\alpha- 1 + (2-\alpha)\varepsilon}$; thus
\[
\sum_{k=1}^\infty k^{-1/\alpha- \varepsilon} \bigl(V_k 1_{V_k \le
k^{1/\alpha+ \varepsilon}} - \mathbf{E}(V_k ;V_k \le k^{1/\alpha+ \varepsilon})\bigr)
\]
is a.s. convergent. Also $\mathbf{E}(V_k; V_k > k^{1/\alpha+
\varepsilon})$ is of order less than $k^{1/\alpha-1 }$ and
$\mathbf{P}(V_k > k^{1/\alpha+ \varepsilon})$ is of order
$k^{-1-\alpha\varepsilon}$ such that $V_k > k^{1/\alpha+
\varepsilon}$ occurs only finitely often a.s. Thus
\[
\sum_{k=1}^\infty k^{-1/\alpha- \varepsilon} (V_k - \gamma)
\]
is a.s. convergent for all $\varepsilon>0$.

For $\beta> \frac1\alpha$, it follows that the sum $\sum_{k=1}^n
k^{-\beta} (V_k-\gamma)$ is a.s. convergent, which is our claim [then
the term $o(n^{1/\alpha- \beta+ \varepsilon})$ is superfluous].
In the case $\beta\le\frac1\alpha$ by Kronecker's lemma a.s.,
\[
\sum_{k=1}^n k^{-\beta} (V_k-\gamma) = o(n^{1/\alpha- \beta+
\varepsilon}) ,
\]
which again is our claim (now $\eta_n$ is superfluous).
\end{pf}

Next recall that $\varsigma$ denotes a random variable with maximally
skewed stable distribution of index $\alpha$ as in~\eqref{stable}. The
following result can be deduced from a general statement on triangular
arrays of independent random variables; see~\cite{fe}, Chapter~XVII,
Section 7; however, a direct proof seems easier.

\begin{lemma} \label{le52}
Let $V_1,V_2,\ldots$ be independent copies of the random variable~\eqref
{Z1}. Then the following holds true:
\begin{longlist}[(ii)]
\item[(i)] Let $1<\alpha< \frac12 (1+ \sqrt5)$. Then
\[
n^{\alpha- 1 - 1/\alpha}\sum_{k=1}^n k^{1-\alpha}(V_k - \gamma)
\stackrel{d}{\to} -c \varsigma,
\]
where
\[
c= \bigl( ( 1+ \alpha-\alpha^2 )\Gamma(2-\alpha)\bigr)^{ -1/\alpha}
.
\]
\item[(ii)]
For $\alpha=\frac12 (1+ \sqrt5)$
\[
(\log n)^{-1/{\alpha}} \sum_{k=1}^n k^{1-\alpha}(V_k - \gamma)
\stackrel{d}{\to} \frac{-\varsigma}{\Gamma(2-\alpha)^{1/\alpha}}
.
\]
\end{longlist}
\end{lemma}

\begin{pf}
(i):
From~\eqref{Z2},~\eqref{Z3} and the theory of stable laws, it follows that
\[
n^{-1/\alpha} (V_1+ \cdots+ V_n-\gamma n ) \stackrel{d}{\to} \frac
{-\varsigma}{\Gamma(2-\alpha)^{1/\alpha}} .\vadjust{\goodbreak}
\]
We express this relation by means of the characteristic functions
$\varphi(u)$ and $e^{\psi(u)}$ of $V-\gamma$ and~$-\varsigma/\Gamma
(2-\alpha)^{1/\alpha}\dvtx\varphi( n^{-1/\alpha}u)^n \to e^{\psi
(u)}$ for all $u \in\mathbb R$, or slightly more generally,
\[
\varphi( v_n n^{-1/\alpha}u)^n \to e^{\psi(u)} ,
\]
if $v_n\to1$. Since $\varphi_n(u)=\varphi(v_n n^{-1/\alpha}u)$ is
again a characteristic function, it follows from Feller~\cite{fe},
Chapter XVII.1, Theorem 1, that for $n \to\infty$
\[
n\bigl ( \varphi(v_nn^{-1/\alpha}u) -1 \bigr) \to\psi(u)
\]
or
\[
\varphi(su) -1 \sim s^{\alpha}\psi(u)\qquad  \mbox{as } s \to0
\]
for all real $u$. Since $\alpha- \alpha^2 > -1$ for $\alpha< \frac
12 (1+ \sqrt5)$, it follows that with $\zeta=(1+ \alpha-\alpha
^2)^{1/\alpha}$
\[
 \sum_{k=1}^n \biggl( \varphi\biggl( \frac{\zeta k^{1-\alpha
}}{n^{1-\alpha+ 1/\alpha}}u\biggr)-1\biggr)\sim\psi(u) \sum
_{k=1}^n \biggl( \frac{\zeta k^{1-\alpha}}{n^{1-\alpha+ 1/\alpha}}
\biggr)^\alpha
\to\psi(u) .
\]
Similarly,
\[
\sum_{k=1}^n \biggl| \varphi\biggl( \frac{\zeta k^{1-\alpha}}{n^{1-\alpha+
1/\alpha}}u\biggr)-1\biggr| \to|\psi(u)|,
\]
and consequently,
\begin{eqnarray*}
&&\sum_{k=1}^n\biggl | \varphi\biggl( \frac{\zeta k^{1-\alpha}}{n^{1-\alpha+
1/\alpha}}u\biggr)-1\biggr|^2 \\
&&\qquad\le\max_{k=1,\ldots,n} \biggl| \varphi\biggl( \frac{\zeta k^{1-\alpha
}}{n^{1-\alpha+ 1/\alpha}}u\biggr)-1\biggr|\sum_{k=1}^n \biggl| \varphi
\biggl( \frac{\zeta k^{1-\alpha}}{n^{1-\alpha+1/\alpha}}u\biggr)-1
\biggr| \to0
\end{eqnarray*}
for $n \to\infty$.

In order to transfer these limit results to characteristic functions we
use that for all complex numbers $z$ with $|z|\le1$,
\[
| z - e^{z-1} | \le c|z-1|^2
\]
for some $c>0$. Therefore, if $|z_1|, \ldots,|z_n| \le1$,
\[
\bigl|z_1 \cdots z_n - e^{(z_1-1)+\cdots+(z_n-1)}\bigr|
\le\sum_{k=1}^n |z_k- e^{z_k-1}|
\le c \sum_{k=1}^n |z_k-1|^2.
\]
We put $z_k=z_{kn}(u)= \varphi( \frac{\zeta k^{1-\alpha
}}{n^{1-\alpha+ 1/\alpha}}u)$. Then the right-hand side goes
to zero, and we obtain
\[
z_{1n}(u) \cdots z_{nn}(u) \to e^{\psi(u) } .
\]
Since the product on the left-hand side is the characteristic function
of $\zeta n^{\alpha- 1 - 1/\alpha}\hspace*{-0.6pt}\times \sum_{k=1}^n k^{1-\alpha} (V_k
- \frac1{\alpha-1} ) $ the claim follows.\vspace*{1pt}\vadjust{\goodbreak}

(ii): This proof goes along the same lines using
\[
\sum_{k=1}^n\biggl( \varphi\biggl( \frac{k^{1-\alpha}}{(\log n)^{1/\alpha}} u\biggr) -1\biggr) \sim\psi(u) \sum_{k=1}^n \biggl( \frac
{k^{1-\alpha}}{(\log n)^{1 /\alpha}}\biggr)^\alpha.
\]
Now $\alpha-\alpha^2 = -1$, thus
\[
\sum_{k=1}^n \biggl(\varphi\biggl( \frac{k^{1-\alpha}}{(\log n)^{1/\alpha}} u\biggr) -1\biggr) \sim\psi(u)\frac1{\log n} \sum_{k=1}^n \frac
1k \sim\psi(u) ,
\]
and the claim follows.
\end{pf}

\section{\texorpdfstring{Proof of Theorem \protect\ref{mainresult}}{Proof of Theorem 1}}\label{sec6}

Again let $2 \le R_1< R_2 < \cdots$ be the points of the stationary
point process $\nu$, and denote
\[
V_j=R_{j+1}-R_j , \qquad j \ge1 .
\]
The random variables $V_1, V_2, \ldots$ are i.i.d. with distribution
\eqref{Z1}.
\begin{lemma} \label{le61}
We have
\[
\int_{[2,n]} x^{1-\alpha} \nu(dx) = \frac{n^{2-\alpha}}{\gamma
(2-\alpha)} - \gamma^{-\alpha}\sum_{k \le n/\gamma} k^{1-\alpha}
(V_k-\gamma) + \delta_n
\]
with
\[
\delta_n= \eta_n+ o_P(n^{1/{\alpha^2} +1-\alpha+ \varepsilon})
\]
for any $\varepsilon>0$, where the random variables $\eta_n$ are
convergent in probability.
\end{lemma}

\begin{pf}
Our starting point is
\[
\int_{[2,n]} x^{1-\alpha} \nu(dx) = \sum_{i=1}^{r_n} R_i^{1-\alpha}
,
\]
where $r_n$ is such that $R_{r_n} \le n < R_{r_n+1}$.
From Lemma~\ref{le51} we have $R_n -\gamma n=o(n^{1/\alpha+
\varepsilon})$ a.s., which implies $r_n-\frac n \gamma=o(n^{
1/\alpha+ \varepsilon})$ a.s.

By a Taylor expansion,
\begin{eqnarray} \label{nu1}
R_i^{1-\alpha} &=& (\gamma i)^{1-\alpha} + (1-\alpha
)(\gamma i)^{-\alpha}(R_i-\gamma i) + \delta_i'
\nonumber
\\[-8pt]
\\[-8pt]
\nonumber
&=& (\gamma i)^{1-\alpha} + (1-\alpha)(\gamma i)^{-\alpha}\sum
_{j=1}^{i-1}(V_j-\gamma) + \delta_i'',
\end{eqnarray}
where the remainder is a.s. of the order
\[
\delta_i'' = O\bigl( i^{-\alpha-1} (R_i- \gamma i)^2\bigr) + O(i^{-\alpha})=
o(i^{2/\alpha- \alpha- 1 + \varepsilon}) .
\]
We consider now the sums of the different terms in~\eqref{nu1}.
\begin{equation} \sum_{i=1}^{r_n} (\gamma i)^{1-\alpha} = \frac{\gamma
^{1-\alpha}}{2-\alpha} r_n^{2-\alpha} + \eta_n' ,\label{nu2}\vadjust{\goodbreak}
\end{equation}
where $\eta_n'$ is a.s. convergent. Further, putting $a_n = (\alpha-1)
\sum_{i>n} i^{-\alpha}$,
\begin{eqnarray*}
(1-\alpha)\sum_{i=1}^{r_n} i^{-\alpha}\sum_{j=1}^{i-1}(V_j-\gamma) &=&
(1-\alpha)\sum_{j=1}^{r_n-1}(V_j-\gamma)\sum_{i=j+1}^{r_n} i^{-\alpha}
\\
&=& a_{r_n}(R_{r_n+1}-R_1-\gamma r_n) - \sum_{j=1}^{r_n}a_j(V_j-\gamma)
.
\end{eqnarray*}
The distribution of $R_{r_n+1}-n$ does not depend on $n$ because of
stationarity, thus $a_{r_n}(R_{r_n+1}-R_1-n)= O_P(n^{1-\alpha})$. Also
$\sum_{j=1}^{n} (a_j-j^{1-\alpha})(V_j-\gamma)$ is a.s. convergent for
$\alpha>1$, since $a_n-n^{1-\alpha}= O(n^{-\alpha})$ and since $V$ has
finite expectation.
It follows
\begin{eqnarray}\label{nu3}
&&(1-\alpha)\sum_{i=1}^{r_n} i^{-\alpha}\sum
_{j=1}^{i-1}(V_j-\gamma)
\nonumber
\\[-8pt]
\\[-8pt]
\nonumber
 &&\qquad= r_n^{1-\alpha}(n-\gamma r_n) -
\sum_{j=1}^{r_n}j^{1-\alpha}(V_j-\gamma)+ \eta_n''+O_P(n^{1-\alpha}) ,
\end{eqnarray}
where $\eta_n''$ is a.s. convergent. Next
\begin{equation} \sum_{i=1}^{r_n} \delta_i'' = \eta_n''' + o(n^{
2/\alpha- \alpha+ \varepsilon})\qquad \mbox{a.s.} \label{nu4}
\end{equation}
for all $\varepsilon>0$, where $\eta_n'''$ is a.s. convergent. Note
that this formula covers two cases: If $\frac2\alpha< \alpha$, then
the sum is a.s. convergent and the right-hand term is superfluous.
Otherwise the term $\eta_n'''$ can be neglected.

Furthermore another Taylor expansion gives
\begin{equation}\qquad
\frac{n^{2-\alpha}}{2-\alpha} = \frac{(\gamma
r_n)^{2-\alpha}}{2-\alpha} + (\gamma r_n)^{1-\alpha}(n-\gamma r_n)+
o(n^{2/\alpha- \alpha+ \varepsilon}) \qquad\mbox{a.s.} \label{nu5}
\end{equation}
Combining~\eqref{nu1} to~\eqref{nu5} gives
\begin{eqnarray}\label{nu6}
\sum_{i=1}^{r_n} R_i^{1-\alpha}&=& \frac{n^{2-\alpha}}{\gamma(2-\alpha)}
- \gamma^{-\alpha} \sum_{j=1}^{r_n} j^{1-\alpha}(V_j-\gamma)
\nonumber
\\[-8pt]
\\[-8pt]
\nonumber
&&{} + \eta_n +
o(n^{2/\alpha-\alpha+ \varepsilon})\qquad
\mbox{a.s.},
\end{eqnarray}
where $\eta_n$ is convergent in probability.

Finally we consider the (loosely notated) difference
\[
\sum_{j=r_n+1}^{n/\gamma} j^{1-\alpha}(V_j-\gamma)=\sum_{j\le n/\gamma}
j^{1-\alpha}(V_j-\gamma)-\sum_{j=1}^{r_n} j^{1-\alpha}(V_j-\gamma) .
\]
For any random sequence of natural numbers $s_n$ such that
$s_n=o(n^{1/\alpha+ \varepsilon})$ a.s. for all $\varepsilon>0$,
\[
\sum_{i \le s_n} (V_i-\gamma) =R_{s_n+1}-R_1-\gamma s_n = o(s_n^{
1/\alpha+ \varepsilon}) =o(n^{1/{\alpha^2}+2\varepsilon+\varepsilon
^2}) \qquad\mbox{a.s.}
\]
Since $r_n-n/\gamma=o(n^{1/\alpha+ \varepsilon})$ a.s. for any
$\varepsilon>0$, this implies for any $\varepsilon>0$ in probability,
\[
\sum_{j=r_n+1}^{n/\gamma}(V_j-\gamma) = o_P(n^{1/{\alpha
^2}+\varepsilon}) .
\]
This implies $\sum_{j=r_n+1}^{n/\gamma} (V_j+\gamma) =o_P(n^{
1/\alpha+ \varepsilon})$. Therefore
\begin{eqnarray*}
&&\Biggl| \sum_{j=r_n+1}^{n/\gamma} j^{1-\alpha}(V_j-\gamma) \Biggr| \\
&&\qquad\le
r_n^{1-\alpha}\Biggl| \sum_{j=r_n+1}^{n/\gamma}(V_j-\gamma) \Biggr| +
\biggl|\biggl(\frac n\gamma\biggr)^{1-\alpha}-r_n^{1-\alpha} \biggr|\sum
_{j=r_n+1}^{n/\gamma} (V_j+\gamma) \\
&&\qquad= o_P(n^{1/{\alpha^2}+1-\alpha+\varepsilon}) + O\bigl(n^{-\alpha}(n-
\gamma r_n)\bigr) o_P(n^{1/\alpha+ \varepsilon})\\
&&\qquad= o_P(n^{1/{\alpha^2}+1-\alpha+\varepsilon}) + o_P(n^{2/\alpha-\alpha+ 2\varepsilon}) .
\end{eqnarray*}
Since $\frac1{\alpha^2} +1 \ge\frac2\alpha$, we end up with
\[
\sum_{j=r_n+1}^{n/\gamma} j^{1-\alpha}(V_j-\gamma) = o_P(n^{1/{\alpha^2} +1-\alpha+ \varepsilon}) .
\]
Combining this estimate with~\eqref{nu6} gives the claim.
\end{pf}

\begin{pf*}{Proof of Theorem \protect\ref{mainresult}}
The total length~\eqref{length} of the $n$-coalescent can be rewritten as
\[
L_n = \sum_{i=0}^{\tau_n-1} \frac{X_{i}}{\rho_{X_i}} E_i,
\]
where $E_0,E_1,\ldots$ denote exponential random variables with
expectation 1, independent among themselves and from the $X_i$.

From Lemma 2.2 in Delmas et al.~\cite{de}, we have for $m \to\infty$
\begin{equation} \rho_m= \frac1{\alpha\Gamma(\alpha)} m^\alpha+
O(m^{\alpha-1}) .
\label{rho}
\end{equation}

In the first step we replace the points $n=X_0> X_1> \cdots$ of a
$\operatorname{CPP}(n)$ by points of a $\operatorname{CPP}(\infty)$: If we take independent versions
of both then for given $\varepsilon> 0$ by Lemma~\ref{le2}, there is a
natural number $s\ge1$ such that with probability at least
$1-\varepsilon$ they meet before $n-s$. From this moment both CPPs can
be coupled. Thus, letting $n \ge X_0'>X_1' >\cdots$ be the points of
the coupled $\operatorname{CPP}(\infty)$ within $[2, n]$, independent of
$E_0,E_1,\ldots,$ and
\[
L_n'= \sum_{i=0}^{\tau_n'-1} \frac{X_{i}'}{\rho_{X_i'}} E_i ,
\]
then due to the the coupling and~\eqref{rho} for $n$ sufficiently big
\[
\mathbf{P}\bigl( |L_n-L_n'| > 3\alpha\Gamma(\alpha) n^{1-\alpha}(E_0+
\cdots+ E_s) \bigr) \le\varepsilon.
\]
Since $\alpha>1$, $L_n-L_n'= o_P(1)$, thus we may replace $L_n$ by
$L_n'$ in our asymptotic considerations.

Thus we work now with a $\operatorname{CPP}(\infty)$ $\mu$, which we couple to a
stationary point process $\nu$ according to~\eqref{bigcoupl6} and \eqref
{bigcoupl7}. Also let $E_0,E_1, \ldots$ be independent of the whole
coupling. We use the formula
\begin{equation}L_n'= \int_{[2,n]} \frac{xE_x}{\rho_x} \mu(dx) ,
\label{Ln}
\end{equation}
in which the exponential random variables now are ordered differently.
Since $\sum_{x \ge1} x^{-\alpha} E_x < \infty$ a.s., it follows from
\eqref{rho} that
\[
L_n' = \alpha\Gamma(\alpha) \int_{[2,n]} \frac{E_x }{x^{\alpha-1}}
\mu(dx) + \eta_{1,n} ,
\]
where $\eta_{1,n}$ is a.s. convergent.

Next $\sum_{x \ge2} x^{-1/2- \varepsilon}(E_x-1)$ is a.s. convergent
for any $\varepsilon>0$. It follows that $\sum_{x \ge2} x^{1-\alpha
}(E_x-1)$ is a.s. convergent for $\alpha> \frac32$ and else a.s. of
order $O(n^{3/2 -\alpha+ \varepsilon})$. Given $\mu$, the same
holds true for
$ \int_{[2,n]} \frac{E_x-1 }{x^{\alpha-1}} \mu(dx)$,
thus
\[
L_n' = \alpha\Gamma(\alpha)\int_{[2,n]} x^{1-\alpha} \mu(dx) + \eta
_{2,n} + o(n^{3/2-\alpha+\varepsilon})\qquad \mbox{a.s.} ,
\]
where again $\eta_{2,n}$ is a.s. convergent.

Next from~\eqref{bigcoupl7} with $2^s < n \le2^{s+1}$
\begin{eqnarray*}
\int_{[2,n]} x^{1-\alpha} \mu(dx) &=& \int_{[2,n]} x^{1-\alpha} \nu
(dx) \\
&&{} + \sum_{r=1}^s \int_{[2^{r-1}+1,2^r]}x^{1-\alpha}
\bigl(\phi_1^{r,\infty}(dx)- \phi_2^{r,\infty}(dx)\bigr)
\\
&&{} + \int_{[2^{s}+1,n]} x^{1-\alpha} \bigl(\phi_1^{s,\infty
}(dx)- \phi_2^{s,\infty}(dx)\bigr).
\end{eqnarray*}
From~\eqref{bigcoupl3} and~\eqref{bigcoupl5} we see that
\begin{eqnarray*}
&&
\biggl|\int_{[2^{r-1}+1,2^r]}x^{1-\alpha} \bigl(\phi
_1^{r,\infty}(dx)- \phi_2^{r,\infty}(dx)\bigr)\biggr| \\
&&\qquad\le2^{r-1} (\alpha
-1) (2^{r-1})^{-\alpha}D_r+ 2 (2^{r-1})^{1-\alpha}D_r ,
\end{eqnarray*}
and the same estimate holds for the last term above. In view of Lemma
\ref{bigcoupl4} and the Borel--Cantelli lemma, we conclude that
\[
\int_{[2,n]} x^{1-\alpha} \mu(dx) = \int_{[2,n]} x^{1-\alpha} \nu
(dx) + \eta_{3,n}
\]
with $\eta_{3,n}$ a.s. convergent. Altogether
\[
L_n' = \alpha\Gamma(\alpha)\int_{[2,n]} x^{1-\alpha} \nu(dx) + \eta
_{4,n} + o(n^{3/2-\alpha+\varepsilon}) ,
\]
where $\eta_{4,n}$ is a.s. convergent. Finally Lemma~\ref{le61} gives a.s.
\begin{eqnarray}\label{exp}
L_n' &=& \frac{\Gamma(\alpha)\alpha(\alpha- 1)}{(2-\alpha
)} n^{2-\alpha} - \Gamma(\alpha)\alpha(\alpha-1)^\alpha\sum_{k \le
 n/\gamma} k^{1-\alpha} (V_k-\gamma)
 \nonumber
 \\[-8pt]
 \\[-8pt]
 \nonumber
 &&{} + \eta
_n+ o_P(n^{1/{\alpha^2} +1-\alpha+ \varepsilon}) + o(n^{3/2-\alpha
+\varepsilon})
\end{eqnarray}
for all $\varepsilon>0$, where $\eta_n$ now is convergent in probability.

We are ready to treat the different cases of Theorem~\ref{mainresult}:

If $1<\alpha< (1+\sqrt5)/2$, then we use that $1/\alpha> 1/\alpha^2$
and $1/\alpha> 1/2$. Therefore the three remainder terms in \eqref
{exp} are all of order $o_P(n^{1/\alpha+1-\alpha})$ and thus may
be neglected. The result follows from an application of Lemma \ref
{le52}. The case $ \alpha= (1+\sqrt5)/2$ is treated in the same way.

If $\alpha> (1+\sqrt5)/2$, then $\frac1{\alpha^2} +1-\alpha< 0$ and
$3/2-\alpha< 0$. Also from Lemma~\ref{le51} it follows that $\sum_{k
\le n/\gamma} k^{1-\alpha} (V_k-\gamma)$ is a.s. convergent.\vspace*{1pt} Thus
it follows from~\eqref{exp} that $L_n' - \frac{\Gamma(\alpha)\alpha
(\alpha- 1)}{(2-\alpha)} n^{2-\alpha}$ is convergent in probability.
To see that the limit of $L_n'$ (and thus $L_n$) is nondegenerate, we
go back to~\eqref{Ln}, respectively,
\begin{eqnarray*}
&& L_n'-\frac{\Gamma(\alpha)\alpha(\alpha- 1)}{(2-\alpha)} n^{2-\alpha
} \\
&&\qquad= \frac{2\alpha\Gamma(\alpha)}{\rho_2} \mu(\{2\}) E_2 +
\biggl(\alpha\Gamma(\alpha)\int_{[3,n]} \frac{xE_x}{\rho_x} \mu(dx)
-\frac{\Gamma(\alpha)\alpha(\alpha- 1)}{(2-\alpha)} n^{2-\alpha} \biggr).
\end{eqnarray*}
As shown the term in brackets in convergent in probability. Also $\mu(\{
2\})=1$ with positive probability. Since the exponential variable $E_2$
is independent from the rest on the right-hand side, the whole limit
has to be nondegenerate. This finishes the proof.
\end{pf*}

%

\printaddresses

\end{document}